# STATISTICAL ESTIMATION IN THE PROPORTIONAL HAZARDS MODEL WITH RISK SET SAMPLING[1]

By Kani Chen

*Hong Kong University of Science and Technology*

Thomas' partial likelihood estimator of regression parameters is widely used in the analysis of nested case-control data with Cox's model. This paper proposes a new estimator of the regression parameters, which is consistent and asymptotically normal. Its asymptotic variance is smaller than that of Thomas' estimator away from the null. Unlike some other existing estimators, the proposed estimator does not rely on any more data than strictly necessary for Thomas' estimator and is easily computable from a closed form estimating equation with a unique solution. The variance estimation is obtained as minus the inverse of the derivative of the estimating function and therefore the inference is easily available. A numerical example is provided in support of the theory.

**1. Introduction.** Thomas' partial likelihood estimate [Thomas (1977) and Oakes (1981)] is the most popular estimate of regression parameters in nested case-control (n-c-c) studies using Cox's proportional hazards model. The partial likelihood score has a simple closed form expression and therefore the estimate is computationally simple with easily available inference. More important, Thomas' estimate only relies on the time-restricted n-c-c data: the failure times of all cases and the covariates of controls (cases) at the time when they are sampled (fail). The aim of this paper is to propose a new estimate that uses only the time-restricted n-c-c data and is more accurate than Thomas' estimate away from the null. This estimate is also easy to compute with readily available inference. Throughout the paper, an estimate is said to be more accurate or efficient than another estimate if the former has smaller asymptotic variance.

Received August 2000; revised August 2003.

[1]Supported in part by an RGC Grant HKUST DAG99/00.SC20.

Statistical analysis of n-c-c designs has attracted considerable attention in the past decade; see Langholz and Thomas (1990, 1991), Goldstein and Langholz (1992), Robins, Rotnitzky and Zhao (1994), Borgan, Goldstein and Langholz (1995), Langholz and Goldstein (1996), Breslow (1996), Samuelsen (1997), Suissa, Edwardes and Biovin (1998), Borgan and Olsen (1999) and Chen (2001), among many others. Some competing estimates were also studied in the literature; see Robins, Rotnitzky and Zhao (1994), Samuelsen (1997) and Chen (2001). However, all these studies depend on the extended n-c-c data which are more than strictly necessary for Thomas' estimate. The extended n-c-c data is defined as the observed failure or censoring times and failure or censoring indices for all cohort members and the entire covariate histories for all cases and controls. Robins, Rotnitzky and Zhao (1994) pointed out that Thomas' estimator is not semiparametrically efficient based on the extended n-c-c data but they only dealt with time-fixed covariates. Samuelsen (1997) proposed an estimator via the inclusion probability method but it is not always more accurate than Thomas' estimator; see discussion in Chen (2001). The estimation method of Chen (2001) leads to a semiparametrically efficient estimator but it inevitably involves estimating and inverting a Fredholm operator, which is computationally difficult.

The most serious practical limitation of the estimators of Samuelsen (1997), Chen (2001) and Robins, Rotnitzky and Zhao (1994), in contrast with Thomas' estimator, is that they rely on the extended n-c-c data rather than only the time-restricted n-c-c data. The extended n-c-c data contain components that are often not available or, even if available, are much less reliable than the time-restricted n-c-c data. First, the nonfailures that are not sampled as controls are usually not closely followed up, and therefore their censoring times are often not accurately observed. This happens particularly when the cohort is loosely defined [Chen and Lo (1999)]. Second, the ascertainment of the entire covariate histories for cases and controls is often too difficult a task to accomplish with reasonable accuracy. Thus, a new estimator would be greatly desirable if it uses only the time-restricted n-c-c data and is reasonably accurate.

The next section introduces notation and Thomas' estimator based on the time-restricted n-c-c data. Section 3 presents the proposed estimator and its consistency and asymptotic normality. A numerical example is provided in Section 4. A few closing remarks are given in Section 5. Proofs are presented in the Appendix.

**2. Thomas' partial likelihood estimator.** Let $\{T, C, Z(\cdot)\}$ be the random triplet of life time, censoring time and covariate process of dimension $d$. Let $Y = \min(T, C)$, $\delta = I(T \leq C)$, $N(t) = \delta I(T \leq t)$ and $Y(t) = I(Y \geq t)$, where $I$ is the indicator function throughout. Consider a cohort of size

$n$. Let $T_i, C_i, Z_i(\cdot), Y_i, \delta_i, N_i(\cdot), Y_i(\cdot), i = 1, 2, \ldots, n$, be the i.i.d. sample analogues. The full cohort data refer to $[(Y_i, \delta_i), \{Z_i(t): t \in [0, Y_i]\}: i = 1, \ldots, n]$. An n-c-c design takes a random sample of size $m$ (for covariate ascertainment) from the risk set at every failure time, excluding the failed subject itself. Let $R_t^*$ denote the index set of a size $m$ random sample selected from all subjects with the minimum of failure and censoring times greater than $t$. The extended n-c-c data refer to $\{(Y_i, \delta_i): i = 1, \ldots, n\} \cup [\{Z_i(t): t \in [0, Y_i]\}, \{Z_j(t): t \in [0, Y_j]\}: \delta_i = 1, j \in R_{Y_i}^*, i, j = 1, \ldots, n]$. The time-restricted n-c-c data refer to $[\{Y_i, Z_i(Y_i), Z_j(Y_i)\}: \delta_i = 1, j \in R_{Y_i}^*, i, j = 1, \ldots, n]$. With the time-restricted n-c-c data, the exact censoring times for the nonfailures are not necessarily specified and the covariate histories for controls (cases) are not observed except at the time when they are sampled (fail).

The Cox proportional hazards model assumes that the conditional hazard of $T$ given $Z$ satisfies

$$\lambda_T(t|Z = z) = \lambda_0(t) \exp\{\beta' z(t)\},$$

where $\beta$ is the parameter to be estimated and $\lambda_0(\cdot)$ is the baseline hazard function. The life time $T$ and the censoring time $C$ are always assumed conditionally independent given $Z$, which is assumed pathwise left continuous with right limit. Cox's partial likelihood estimator of $\beta$ based on the full cohort data, denoted by $\hat{\beta}_C$, is the solution of

$$U_C(\beta) \equiv \sum_{i=1}^n \int_0^\tau \left[ Z_i(t) - \frac{\sum_{j \in R_t} Z_j(t) \exp\{\beta' Z_j(t)\}}{\sum_{j \in R_t} \exp\{\beta' Z_j(t)\}} \right] dN_i(t) = 0,$$

where $\tau = \sup\{t: \text{pr}(Y > t) > 0\}$ and $R_t$ is the risk set at time $t$; that is, $R_t = \{j: Y_j \geq t, j = 1, \ldots, n\}$. Thomas' partial likelihood estimator of $\beta$, denoted by $\hat{\beta}_P$, is the solution of

$$U_P(\beta) \equiv \sum_{i=1}^n \int_0^\tau \left[ Z_i(t) - \frac{\sum_{j \in R_t^* \cup \{i\}} Z_j(t) \exp\{\beta' Z_j(t)\}}{\sum_{j \in R_t^* \cup \{i\}} \exp\{\beta' Z_j(t)\}} \right] dN_i(t) = 0.$$

It is clear that Thomas' estimator uses only the time-restricted n-c-c data. It is proved in Goldstein and Langholz (1992) (where $\tau$ is set to be 1 for convenience) that, under certain regularity conditions,

$$n^{1/2}(\hat{\beta}_P - \beta) \to N(0, \Sigma_P^{-1}),$$

where $\Sigma_P = \Sigma_C - \Sigma_a$,

$$\Sigma_C = \int_0^\tau E[\{Z(t) - \mu(t)\}^{\otimes 2} \exp\{\beta' Z(t)\} Y(t)] \lambda_0(t) \, dt,$$

$$\Sigma_a = \frac{1}{m+1} \int_0^\tau E\left( \frac{[\sum_{j=1}^{m+1} \{Z_j(t) - \mu(t)\} \exp\{\beta' Z_j(t)\}]^{\otimes 2}}{\sum_{j=1}^{m+1} \exp\{\beta' Z_j(t)\}} \right|$$

$$\left. Y_1 \geq t, \ldots, Y_{m+1} \geq t \right) \text{pr}(Y \geq t) \lambda_0(t) \, dt$$

and

$$\mu(t) = E\{Z(t)|Y = t, \delta = 1\} = \frac{E[Z(t)\exp\{\beta'Z(t)\}Y(t)]}{E[\exp\{\beta'Z(t)\}Y(t)]}.$$

Here and throughout the paper, $v^{\otimes 2} = vv'$ for any vector $v$ of dimension $d$. Moreover, $-\dot{U}_P(v)|_{v=\hat{\beta}}/n$ is a consistent estimator of $\Sigma_P$. It is also well known that $n^{1/2}(\hat{\beta}_C - \beta) \to N(0, \Sigma_C^{-1})$. Throughout the paper, the true value of $\beta$ is still denoted by $\beta$.

REMARK. Set

$$\Psi(t) = \sum_{j=1}^{m+1}(\{Z_j(t) - \mu(t)\}\exp\{\beta'Z_j(t)\}) \Big/ \sum_{j=1}^{m+1}\exp\{\beta'Z_j(t)\}$$

for ease of notation. It follows from the expression of $\Sigma_P$ given in Goldstein and Langholz (1992) that

$$\Sigma_P = \int_0^\tau E([Z_1(t) - \mu(t) - \Psi(t)]^{\otimes 2}$$

$$\times \exp\{\beta'Z_1(t)\}Y_1(t)|Y_2 \geq t, \ldots, Y_{m+1} \geq t)\lambda_0(t)\,dt$$

$$= \int_0^\tau E([Z_1(t) - \mu(t) - \Psi(t)]^{\otimes 2}$$

$$\times \exp\{\beta'Z_1(t)\}|Y_1 \geq t, \ldots, Y_{m+1} \geq t)\operatorname{pr}(Y \geq t)\lambda_0(t)\,dt$$

$$= \int_0^\tau E([Z_1(t) - \mu(t)]^{\otimes 2}\exp\{\beta'Z_1(t)\}Y_1(t))\lambda_0(t)\,dt$$

$$- \int_0^\tau E([Z_1(t) - \mu(t)]\Psi(t)'$$

$$\times \exp\{\beta'Z_1(t)\}|Y_1 \geq t, \ldots, Y_{m+1} \geq t)\operatorname{pr}(Y \geq t)\lambda_0(t)\,dt$$

$$- \int_0^\tau E(\Psi(t)[Z_1(t) - \mu(t)]'$$

$$\times \exp\{\beta'Z_1(t)\}|Y_1 \geq t, \ldots, Y_{m+1} \geq t)\operatorname{pr}(Y \geq t)\lambda_0(t)\,dt$$

$$+ \int_0^\tau E(\Psi(t)^{\otimes 2}\exp\{\beta'Z_1(t)\}|Y_1 \geq t, \ldots, Y_{m+1} \geq t)\operatorname{pr}(Y \geq t)\lambda_0(t)\,dt$$

$$= \Sigma_C - \Sigma_a - \Sigma_a + \Sigma_a = \Sigma_C - \Sigma_a.$$

The expression of $\Sigma_P$ in the first line is intuitively well understandable from the partial likelihood nature of the Thomas estimator $\hat{\beta}_P$.

**3. A new estimator and its inference.** The motivation of the new estimator is described in the following. Observe that the ratio in the expression of $U_P$ can be viewed as an estimator of $\mu(t)$. This estimator, although unbiased, uses only $m+1$ observations: $m$ controls plus 1 case. Heuristically, its estimation and accuracy can be improved by utilizing more observations of relevance. One way to do so is to consider altogether the controls in $R_s^*$ for all failure times $s$ in a neighborhood of $t$. With a more accurate estimation and approximation of $\mu(t)$ plus a proper weighting scheme, one can presumably construct a new estimator of the regression parameters. The details of the construction are as follows. Set $\bar{N}(s) = (1/n)\sum_{i=1}^n N_i(s)$, $b(t) = E[\exp\{\beta' Z(t)\}|Y \geq t]$,

$$\tilde{w}(t) = \frac{mb(t)}{\exp\{\beta'Z(t)\} + mb(t)}, \qquad g(t) = \frac{E[\tilde{w}(t)Z(t)\exp\{\beta'Z(t)\}Y(t)]}{E[\tilde{w}(t)\exp\{\beta'Z(t)\}Y(t)]},$$

and let $\tilde{w}_i(t)$ be the i.i.d. copies of $\tilde{w}(t)$. Let $\psi_n(x) = \psi(n^{1/3}x)$, $x \in (-\infty, \infty)$, where $\psi$ is an infinitely differentiable nonnegative even function with bounded support. For ease of notation, suppose the support of $\psi$ is $(-1, 1)$. Use

$$(3.1) \qquad \hat{b}(t) \equiv \frac{\int_0^\tau \sum_{j \in R_s^*} \exp\{\hat{\beta}_P' Z_j(s)\} \psi_n(t-s) \, d\bar{N}(s)}{m \int_0^\tau \psi_n(t-s) \, d\bar{N}(s)}$$

to approximate $b(t)$ and

$$(3.2) \qquad w_i(t) \equiv \frac{m\hat{b}(t)}{\exp\{\hat{\beta}_P' Z_i(t)\} + m\hat{b}(t)}$$

to approximate $\tilde{w}_i(t)$, $i = 1, 2, \ldots$. Let

$$S_k(t, \beta) = \int_0^\tau \sum_{j \in R_s^*} w_j(s) Z_j^k(s) \exp\{\beta' Z_j(s)\} \psi_n(t-s) \, d\bar{N}(s), \qquad k = 0, 1, 2.$$

Throughout the paper, the power 2 on covariates $Z$ or $Z_j, j \geq 1$, always means the outer product $\otimes 2$. The notion $A \geq B$ for any two $d \times d$ nonnegative definite matrices means that $A - B$ is nonnegative definite. The proposed estimator, denoted by $\hat{\beta}$, is the solution of

$$(3.3) \qquad U(\beta) \equiv \sum_{i=1}^n \int_0^\tau w_i(t) \left\{ Z_i(t) - \frac{S_1(t, \beta)}{S_0(t, \beta)} \right\} dN_i(t) = 0.$$

We note that $S_1(t, \beta)/S_0(t, \beta)$ may be viewed as an estimator of $g(t)$, a weighted version of $\mu(t)$. Suppose, in the definition of $S_k(t, \beta)$, we use 1 instead of the presently defined $w_i$ as weights. Then $S_1(t, \beta)/S_0(t, \beta)$ is an estimator of $\mu(t)$ which should heuristically be more accurate than its counterpart in $U_P(\beta)$. The present choice of weights is optimal in the sense that no other choices will produce estimators of $\beta$ with smaller asymptotic variance; see further discussion in Section 5. As a result, (3.3) may produce more accurate estimators than Thomas' estimator. The difference between the weights $w_i(t)$ used here and those used in Sasieni (1993b) is that $w_i(t)$ depends on the covariate $Z_i(t)$ while those in Sasieni (1993b) do not.

Denote by $B_c$ the closed ball in $d$ dimensional real space centered at the origin with radius $c > 0$, where $c$ is a large but fixed constant. Let $C^\infty(B_c)$ denote the set of all infinitely differentiable functions defined on $B_c$. Some regularity conditions are assumed here:

(i) $\mathrm{pr}\{\sup_{t \in [0,\tau]} |Z(t)| \leq c\} = 1$;

(ii) the baseline hazard function $\lambda_0(t)$ is bounded away from 0 and infinity on $[0,\tau]$ and has continuous second derivative;

(iii) $\mathrm{pr}\{Y(t) = 1\} > 0$ for all $t \in [0,\tau]$;

(iv) $\Sigma_P$ and $\Sigma$ defined in (3.5) are positive definite;

(v) for any $\phi(\cdot) \in C^\infty(B_c)$, $E[\phi\{Z(t)\}Y(t)]$ as a function of $t$ on $[0,\tau]$ has continuous second derivative;

(vi) for any $\phi(\cdot) \in C^\infty(B_c)$, the process $n^{-1/2} \sum_{i=1}^n (\phi\{Z_i(t)\} Y_i(t) - E[\phi\{Z(t)\}Y(t)])$, as a process of $t$, converges to a Gaussian process on $[0,\tau]$ as $n \to \infty$.

THEOREM. *Assume the above conditions* (i)–(vi) *hold. Then*

$$(3.4) \qquad n^{1/2}(\hat{\beta} - \beta) \to N(0, \Sigma^{-1}),$$

*where*

$$(3.5) \qquad \Sigma = \int_0^\tau E[\tilde{w}(t)\{Z(t) - g(t)\}^{\otimes 2} \exp\{\beta' Z(t)\} Y(t)] \lambda_0(t)\, dt.$$

*Moreover,* $-(1/n)\dot{U}(v)|_{v=\hat{\beta}}$ *is a consistent estimator of* $\Sigma$, *where*

$$\dot{U}(v)|_{v=\hat{\beta}} = -\int_0^\tau \left[ \frac{S_2(t,\hat{\beta})}{S_0(t,\hat{\beta})} - \left\{ \frac{S_1(t,\hat{\beta})}{S_0(t,\hat{\beta})} \right\}^{\otimes 2} \right] \sum_{i=1}^n w_i(t)\, dN_i(t).$$

REMARK. There are six conditions assumed in Goldstein and Langholz (1992) to ensure the consistency and asymptotic normality of $\hat{\beta}_P$. Conditions (i)–(iv) here are analogous to Conditions 2–5 in Goldstein and Langholz (1992). Their Conditions 1 and 6 are implied in the model description in Section 2 and therefore are not listed here. The inheritance of their conditions is understandable since the proposed estimator $\hat{\beta}$ uses the Thomas estimator $\hat{\beta}_P$. Conditions (i)–(iii) and (v) are necessary for using the empirical approximations (e.g., the proof of part 1 of the Lemma in the Appendix) to obtain the rates of convergence for various random quantities, as in the Lemma in the Appendix. Condition (iii) here, parallel to Condition 4 in Goldstein and Langholz (1992), can be relaxed with increasing technicalities involving the tail behavior near the endpoint $\tau$. Condition (iv) validates the asymptotic normality of $\hat{\beta}$ claimed in (3.4). It is also used in proving the consistency of $\hat{\beta}$; see Step 3 of the proof of the Theorem in the Appendix. The requirement of differentiability in conditions (ii) and (v), for obtaining the bounds of kernel estimation, may appear to be restrictive in that it does not allow, for example, $\text{pr}(Y \geq t)$ and $E\{Z(t)\}$ to be discontinuous, although pathwise discontinuous $Z(\cdot)$ are not excluded. In fact, the differentiability requirements in conditions (ii) and (v) can be relaxed to piecewise differentiability. Then, functions such as $\text{pr}(Y \geq t)$ and $E\{Z(t)\}$ can be discontinuous at a finite number of time points. The important case of fixed censoring is also covered. The proof of the Theorem under such relaxed conditions requires a careful but regular treatment on the edge effect caused by the discontinuity points and the endpoints, which is partly reflected in the Lemma. The current presentation was chosen to avoid a lengthy but not essential technical argument. In particular, condition (v) is satisfied when $\int P\{C \geq t | Z(t) = z\} \phi(z) f_t(z) \, d\mathcal{L}(z)$ as a function of $t$ on $[0, \tau]$ is piecewise twice continuously differentiable, where $\phi(\cdot) \in C^\infty(B_c)$ and $f_t(z)$ is the density of $Z(t)$ with respect to a measure $\mathcal{L}$ which can be a combination of the Lebesgue measure and a counting measure. Condition (vi) is essentially about the tightness of the sequence, which can be ensured by a certain Lipschitz condition on the increments $Z(t) - Z(s)$. For example, it is satisfied if there exist an $a > 1$ and $A > 0$ such that $E(|Z(t) - Z(s)|^2) \leq A|t - s|^a$ for all $s, t \in [0, \tau]$. More relaxed but technical conditions in terms of metric entropy may be found in Pollard (1990) or van der Vaart and Wellner (1996). Condition (vi) is used to obtain uniform bounds on $[0, \tau]$ for sequences of random processes in concern; see (A.13). Again, a piecewise version of condition (vi) is sufficient.

REMARK. One can use another version of the weights $w_i(t)$ by replacing $\hat{\beta}_P$ by $\beta$ in the definitions of $\hat{b}$ and $w_i$ in (3.1) and (3.2). The advantage is that one does not have to compute $\hat{\beta}_P$ first to obtain $\hat{\beta}$. The disadvantage is that the estimating equation (3.3) may possibly have multiple roots. Still, the same proof shows that one of the roots is consistent and asymptotically normal with the same asymptotic variance $\Sigma^{-1}$.

PROPOSITION. *Assume the conditions of the Theorem hold. Then $\Sigma_P^{-1} \geq \Sigma^{-1}$ and equality holds if and only if $\beta = 0$.*

This inequality justifies that the asymptotic variance of the proposed estimator is smaller than that of Thomas' estimator away from the null.

**4. A numerical example.** Some simulation results are presented in this section. The covariate process is such that $Z(t) = 4tu_1 + u_2$, $t \in [0, 1]$, where

TABLE 1
*Summary of the simulation results*[a]

|  | $\beta = 0$ | | | $\beta = 1$ | | | $\beta = 2$ | | |
|---|---|---|---|---|---|---|---|---|---|
|  | Ave Est | Emp Var | Est Var | Ave Est | Emp Var | Est Var | Ave Est | Emp Var | Est Var |
| Fixed censoring: $C = 1$ | | | | | | | | | |
| Cen Prop |  | 0.368 |  |  | 0.603 |  |  | 0.667 |  |
| Cox | −0.004 | 0.013 | 0.012 | 1.012 | 0.030 | 0.029 | 2.025 | 0.067 | 0.063 |
| $m = 1$ | | | | | | | | | |
| Thomas | −0.003 | 0.026 | 0.025 | 1.039 | 0.089 | 0.087 | 2.127 | 0.316 | 0.299 |
| Proposed | −0.004 | 0.025 | 0.027 | 0.989 | 0.065 | 0.081 | 1.962 | 0.162 | 0.197 |
| $m = 2$ | | | | | | | | | |
| Thomas | −0.003 | 0.019 | 0.018 | 1.034 | 0.059 | 0.058 | 2.096 | 0.189 | 0.174 |
| Proposed | −0.003 | 0.019 | 0.019 | 1.031 | 0.049 | 0.056 | 2.028 | 0.123 | 0.140 |
| $m = 3$ | | | | | | | | | |
| Thomas | −0.001 | 0.017 | 0.016 | 1.034 | 0.054 | 0.048 | 2.073 | 0.142 | 0.134 |
| Proposed | −0.001 | 0.017 | 0.017 | 1.039 | 0.047 | 0.048 | 2.041 | 0.106 | 0.116 |
| Random censoring: $C\|Z \sim U[0, \min(1, \|Z(0.25)\|)]$ | | | | | | | | | |
| Cen Prop |  | 0.771 |  |  | 0.849 |  |  | 0.843 |  |
| Cox | −0.007 | 0.059 | 0.056 | 1.009 | 0.080 | 0.079 | 2.018 | 0.124 | 0.121 |
| $m = 1$ | | | | | | | | | |
| Thomas | 0.009 | 0.138 | 0.123 | 1.106 | 0.336 | 0.301 | 2.231 | 0.845 | 0.880 |
| Proposed | −0.001 | 0.126 | 0.147 | 0.991 | 0.195 | 0.276 | 1.963 | 0.376 | 0.453 |
| $m = 2$ | | | | | | | | | |
| Thomas | −0.011 | 0.094 | 0.087 | 1.063 | 0.189 | 0.176 | 2.159 | 0.481 | 0.435 |
| Proposed | −0.013 | 0.095 | 0.095 | 1.031 | 0.145 | 0.176 | 2.024 | 0.256 | 0.300 |
| $m = 3$ | | | | | | | | | |
| Thomas | −0.003 | 0.082 | 0.078 | 1.069 | 0.154 | 0.144 | 2.137 | 0.370 | 0.334 |
| Proposed | −0.006 | 0.082 | 0.082 | 1.057 | 0.133 | 0.147 | 2.043 | 0.205 | 0.243 |

[a]"Ave Est" and "Emp Var" stand for the averages and empirical variances of the estimates over 2000 simulations. "Est Var" stands for the average of the estimated variances over 2000 simulations. "Cen Prop" stands for the proportion of censoring. "Cox," "Thomas" and "Proposed" refer, respectively, to the Cox estimate based on full cohort data, Thomas' estimate and the proposed estimates based on time-restricted n-c-c data.

$u_1$ and $u_2$ are two independent random variables uniformly distributed on $[-1, 1]$. The baseline hazard function $\lambda_0(t)$ is set to be constant at 1. We consider separately two different types of censorship: the fixed censoring with $\mathrm{pr}(C = 1) = 1$, and the random censoring with the conditional distribution of $C$ given that $Z$ is the uniform distribution on $[0, \min(1, |Z(0.25)|)]$. The function $\psi_n(t)$ is chosen to be $I(|t| \leq 0.05)$. The regression parameter $\beta$ takes values 0, 1 and 2 and the size of controls to be selected from each risk set is 1, 2 and 3. The sample size is 200. For each scenario, 2000 simulations are conducted and Thomas' estimator and the proposed estimator are calculated. For reference, we also calculate Cox's partial likelihood estimator based on full cohort data. The results are presented in Table 1. Table 1 shows that the proposed estimate has indeed smaller variances (and mean squared errors as well) than Thomas' estimate if $\beta \neq 0$. When $\beta = 0$, the two estimates have about the same asymptotic variances. These simulation results are consistent with the Theorem and the Proposition. In this example, when $\beta = 2$ the bias of Thomas' estimate appears to be relatively serious, while that of the proposed estimate is always negligible. We also notice that, in a few cases of this example, the variance estimation appears to be biased down for Thomas' estimate and biased up for the proposed estimate. Typically, the latter will result in conservative but still valid inferences. As sample size increases, the bias tends to be negligible. It is concluded that this simulation example provides solid evidence in support of the established theoretical results.

**5. Closing remarks.** In summary, the proposed estimate $\hat{\beta}$ is asymptotically more accurate than Thomas' estimate away from the null and it uses only the time-restricted n-c-c data which is strictly necessary for Thomas' estimate. This estimate is relatively easy to compute: The estimating equation takes a simple closed form and has a unique solution. Its inference is equally easy to obtain as the variance estimate is simply minus the inverse of the derivative of the estimating function. Unlike the case of curve estimation, the problem of (optimal) bandwidth choice is much less significant here. In the definition of $\psi_n$, the order of bandwidth is $n^{-1/3}$. In fact, with little modification of the proof, the Theorem holds for all bandwidths of order $n^{-r}$ with $r \in (1/4, 1/2)$. At least within this range, the choice of bandwidth does not affect the first-order asymptotic behavior of the estimate. In practice, however, a proper objective or data-driven choice of bandwidth should be valuable for the implementation of the estimation procedure.

Although it is not clear whether the proposed estimator is semiparametric efficient based on the time-restricted n-c-c data, it does have the following optimality. Consider the class of estimators as solutions of

$$\sum_{i=1}^{n} \int_0^{\tau} \left[ h\{Z_i(t), t, \beta\} - \frac{\int_0^{\tau} \sum_{j \in R_s^*} h\{Z_j(s), s, \beta\} \exp\{\beta' Z_j(s)\} \psi_n(t-s) \, d\bar{N}(s)}{\int_0^{\tau} \sum_{j \in R_s^*} \exp\{\beta' Z_j(s)\} \psi_n(t-s) \, d\bar{N}(s)} \right] dN_i(t) = 0,$$

where $h$ is any bounded infinitely differentiable function. Heuristically, if $h\{Z_j(t), t, \beta\}$ is replaced by $w_j(t)\{Z_j(t) - S_1(t, \beta)/S_0(t, \beta)\}$, it can be shown that the above estimating function is $\sum_{i=1}^{n} \int_0^\tau h\{Z_i(t), t, \beta\} \, dN_i(t) + o_P(n^{1/2})$. Hence the above equation with this particular choice of $h$ is asymptotically equivalent to (3.3) in the sense that the resulting estimators of $\beta$ are asymptotically equivalent. In this asymptotic sense, (3.3) might be viewed as approximately a member of the above general class of estimating equations. More important, it can be proved under regularity conditions that the theoretical optimal choice of $h\{Z_j(t), t, \beta\}$ is $\tilde{w}_j(t)\{Z_j(t) - g(t)\}$, which in actual construction is approximated by $w_j(t)\{Z_j(t) - S_1(t,\beta)/S_0(t,\beta)\}$. It implies that no other choice of $h$ used in the above estimating equation shall result in estimators of $\beta$ with asymptotic variance smaller than $\Sigma^{-1}$ and that no other choice of weights used in (3.3) will produce estimators of $\beta$ with asymptotic variance smaller than $\Sigma^{-1}$.

## APPENDIX

**Proofs of the Theorem and the Proposition.** More notation is needed. Throughout the Appendix, the notion $|\cdot|$ for a vector or matrix means the sum of the absolute values of all elements. Set $n_t = \sum_{i=1}^n I(Y_i > t)$, $f(t) = E[\exp\{\beta' Z(t)\} Y(t)]\lambda_0(t)$ and $h_k(t) = E[\tilde{w}(t) Z^k(t) \exp\{\beta' Z(t)\} | Y \geq t]$, $k = 0, 1, 2$. Notice that together conditions (i), (ii) and (v) ensure that $f(t)$ is bounded above, bounded above 0 and has a continuous second derivative on $[0, \tau]$. Much of the proof relies on counting process martingale techniques; see, for example, Andersen and Gill (1982). The following lemma provides the approximations used in the proof of the Theorem.

LEMMA. *Assume conditions* (i)–(vi) *hold. Let $\varepsilon > 0$ be an arbitrary number.*

1. *Let $f_n(t) = n^{1/3} E\{\psi_n(t - Y)\delta\}$ and $a_n(t) = n^{1/3} E\{\delta(Y - t)\psi_n(Y - t)\}$. Then*

$$\text{(A.1)} \quad \sup_{0 \leq t \leq \tau} \left| n^{1/3} \int_0^\tau \psi_n(t - s) \, d\bar{N}(s) - f_n(t) \right| = O_P(n^{-1/3+\varepsilon}),$$

$$\text{(A.2)} \quad \sup_{0 \leq t \leq \tau} \left| n^{1/3} \int_0^\tau (s - t) \psi_n(t - s) \, d\bar{N}(s) - a_n(t) \right| = O_P(n^{-2/3+\varepsilon}).$$

*Moreover, $f_n(\cdot)$ and $a_n(\cdot)$ are continuous on $[0, \tau]$, satisfying*

$$0 < \inf_n \inf_{t \in [0,\tau]} f_n(t) \leq \sup_n \sup_{t \in [0,\tau]} f_n(t) < \infty,$$

(A.3)

$$\sup_{t \in [0,\tau]} |a_n(t)| = O(n^{-1/3}),$$

*and $f_n(t) = f(t)$ and $a_n(t) = 0$ for $t \in (n^{-1/3}, \tau - n^{-1/3})$.*

2. Let $\tilde{S}_k(t,\beta)$ be defined the same as $S_k(t,\beta)$, $k=0,1,2$, except with $w_i$ replaced by $\tilde{w}_i$. Then

$$\sup_{t\in[0,\tau]} |\hat{b}(t) - b(t)| = O_P(n^{-1/3+\varepsilon}); \tag{A.4}$$

$$\sup_{1\leq i\leq n} \sup_{t\in[0,\tau]} |w_i(t) - \tilde{w}_i(t)| = O_P(n^{-1/3+\varepsilon}); \tag{A.5}$$

$$\sup_{t\in[0,\tau]} \left| \frac{S_k(t,\beta) - \tilde{S}_k(t,\beta)}{m \int_0^\tau \psi_n(t-s)\, d\bar{N}(s)} \right| = O_P(n^{-1/3+\varepsilon}); \tag{A.6}$$

$$\sup_{t\in[0,\tau]} \left| \frac{\tilde{S}_k(t,\beta)}{m \int_0^\tau \psi_n(t-s)\, d\bar{N}(s)} - E[\tilde{w}(t)Z^k(t)\exp\{\beta'Z(t)\}|Y>t] \right| = O_P(n^{-1/3+\varepsilon}); \tag{A.7}$$

$$\sup_{t\in[0,\tau]} \left| \frac{\tilde{S}_1(t,\beta)}{\tilde{S}_0(t,\beta)} - g(t) \right| = O_P(n^{-1/3+\varepsilon}); \tag{A.8}$$

$$\sup_{s\in[0,\tau]} \left| \frac{1}{n} \int_0^\tau \psi_n(t-s) \frac{\sum_{i=1}^n \tilde{w}_i(t)\, dN_i(t)}{\tilde{S}_0(t,\beta)} - \frac{\lambda_0(s)P(Y\geq s)}{mf_n(s)} \right| = O_P(n^{-1/3+\varepsilon}), \tag{A.9}$$

where $\lambda_0(s)P(Y\geq s)/\{mf_n(s)\} = 1/\{mb(s)\}$ for $s\in(n^{-1/3}, \tau - n^{-1/3})$. Moreover, (A.7)–(A.9) also hold if $\tilde{S}_k(t,\beta)$ is replaced by $S_k(t,\beta)$.

PROOF. 1. Observe that $\mathrm{var}\{\psi_n(t-Y)\delta\} \leq \int_0^\tau \psi_n^2(t-s)f(s)\, ds = O(n^{-1/3})$ and write

$$n^{1/3} \left| n^{1/3} \int_0^\tau \psi_n(t-s)\, d\bar{N}(s) - f_n(t) \right|$$
$$= n^{-1/3} \left| \sum_{i=1}^n [\psi_n(t-Y_i)\delta_i - E\{\psi_n(t-Y)\delta\}] \right|.$$

Set $M$ as a large but fixed number. It follows from Bernstein's inequality [see, e.g., van der Vaart and Wellner (1996), page 102] that

$$P\left( n^{-1/3} \left| \sum_{i=1}^n [\psi_n(t-Y_i)\delta_i - E\{\psi_n(t-Y)\delta\}] \right| \geq n^\varepsilon \right) \leq 2n^{-M}$$

for all large $n$. Set $\mathcal{A}_n = \{k\tau/n^2 : k = 0, 1, \ldots, n^2\}$. The above exponential inequality ensures, through the Borel–Cantelli lemma, that

$$\sup_{t \in \mathcal{A}_n} n^{-1/3} \left| \sum_{i=1}^{n} [\psi_n(t - Y_i)\delta_i - E\{\psi_n(t - Y)\delta\}] \right| = O(n^\varepsilon)$$

almost surely. Extending the supremum over $\mathcal{A}_n$ to over $[0, \tau]$, the above equality still holds by the differentiability of the kernel function $\psi$. Therefore (A.1) holds. (A.2) can be proved in a similar fashion. (A.3) and the rest of the claims can be verified by direct calculation using Taylor expansion.

2. Set $d_j(t) = [\exp\{\beta' Z_j(t)\} - b(t)]Y_j(t+)$. Observe (A.1)–(A.3) and that $\hat{\beta}_P - \beta = O_P(n^{-1/2})$. One can apply Taylor expansion and write

$$\hat{b}(t) - b(t)$$

$$= \frac{1}{m \int_0^\tau \psi_n(t - s) \, d\bar{N}(s)}$$

$$\times \left[ \int_0^\tau \sum_{j \in R_s^*} [\exp\{\hat{\beta}_P' Z_j(s)\} - b(s)]\psi_n(t - s) \, d\bar{N}(s) \right.$$

$$\left. \qquad + m \int_0^\tau \{b(s) - b(t)\}\psi_n(t - s) \, d\bar{N}(s) \right] \qquad \text{(A.10)}$$

$$= \frac{1}{m \int_0^\tau \psi_n(t - s) \, d\bar{N}(s)}$$

$$\times \left[ \int_0^\tau \sum_{j \in R_s^*} [\exp\{\hat{\beta}_P' Z_j(s)\} - \exp\{\beta' Z_j(s)\}]Y_j(s+)\psi_n(t - s) \, d\bar{N}(s) \right.$$

$$+ \int_0^\tau \sum_{j \in R_s^*} [\exp\{\beta' Z_j(s)\} - b(s)]Y_j(s+)\psi_n(t - s) \, d\bar{N}(s)$$

$$\left. + m\dot{b}(t) \int_0^\tau (s - t)\psi_n(t - s) \, d\bar{N}(s) \right]$$

$$+ O_P(n^{-2/3})$$

$$= (\hat{\beta}_P - \beta)' E[Z(t) \exp\{\beta' Z(t)\} | Y \geq t]$$

$$+ \frac{1}{m f_n(t)} \int_0^\tau \sum_{j=1}^{n} \left\{ I(j \in R_s^*) - \frac{m}{n_s} \right\} d_j(s) n^{1/3} \psi_n(t - s) \, d\bar{N}(s)$$

$$+ \frac{1}{n f_n(t) P(Y \geq t)} \sum_{j=1}^{n} \int_0^\tau d_j(s) n^{1/3} \psi_n(t - s) f(s) \, ds$$

$$+ \dot{b}(t) a_n(t) / f_n(t) + o_P(n^{-1/2}).$$

Here and in the following, $O_P(\cdot)$ and $o_P(\cdot)$ are uniform over $t \in [0, \tau]$. The first term in the last expression is $O_P(n^{-1/2})$ by the asymptotic normality of $\hat{\beta}_P$ established in Goldstein and Langholz (1992). Let $\mathcal{F}$ denote the $\sigma$-algebra generated by $[\{Y_i, \delta_i, Z_i(\cdot)\}, i = 1, 2, \ldots]$. The integrands of the second term are conditionally independent with conditional mean zero when conditioning on $\mathcal{F}$. Thus, the empirical approximation analogous to the proof of (A.1) can be applied to show the second term is $O_P(n^{-1/3+\varepsilon})$. Since $d_j(t)$ is uniformly bounded with mean zero, the third term can be similarly shown to be $O_P(n^{-1/3+\varepsilon})$. The fourth term is $O_P(n^{-1/3+\varepsilon})$ by part 1. Therefore (A.4) is proved.

To show (A.5), apply the mean value theorem and write

$$(A.11) \quad \begin{aligned} w_i(t) - \tilde{w}_i(t) \\ = \tilde{w}_i(t)\{1 - \tilde{w}_i(t)\}[\{\hat{b}(t) - b(t)\} - (\hat{\beta}_P' - \beta')Z_i(t)] + o_P(n^{-1/2}), \end{aligned}$$

where $o_P(\cdot)$ is uniform over $[0, \tau]$. Then (A.5) follows from (A.4), the boundedness of $Z(\cdot)$ in condition (i) and the asymptotic normality of $\hat{\beta}_P$.

Equation (A.6) follows directly from (A.5) and the definitions of $S_k(t, \beta)$ and $\tilde{S}_k(t, \beta)$.

To show (A.7), let $\hat{h}_k(t) = (1/n_t) \sum_{j=1}^n \tilde{w}_j(t) Z_j^k(t) \exp\{\beta' Z_j(t)\} Y_j(t+)$ and recall the definition of $h_k(\cdot)$. Conditions (v) and (vi) imply that $\sup_{t \in [0,\tau)} |\hat{h}_k(t) - h_k(t)| = O_P(n^{-1/2})$. Therefore we can write

$$\begin{aligned} & \frac{\tilde{S}_k(t, \beta)}{m \int_0^\tau \psi_n(t-s)\, d\bar{N}(s)} - h_k(t) \\ &= \frac{\int_0^\tau \sum_{j \in R_s^*}[\tilde{w}_j(s) Z_j^k(s) \exp\{\beta' Z_j(s)\} - h_k(t)] \psi_n(t-s)\, d\bar{N}(s)}{m \int_0^\tau \psi_n(t-s)\, d\bar{N}(s)} \\ &= \frac{\int_0^\tau \sum_{j \in R_s^*}[\tilde{w}_j(s) Z_j^k(s) \exp\{\beta' Z_j(s)\} - h_k(s)] \psi_n(t-s)\, d\bar{N}(s)}{m \int_0^\tau \psi_n(t-s)\, d\bar{N}(s)} \\ &\quad + O_P(n^{-1/3+\varepsilon}) \\ &= \frac{\int_0^\tau \sum_{j \in R_s^*}[\tilde{w}_j(s) Z_j^k(s) \exp\{\beta' Z_j(s)\} - \hat{h}_k(s)] \psi_n(t-s)\, d\bar{N}(s)}{m \int_0^\tau \psi_n(t-s)\, d\bar{N}(s)} \\ &\quad + O_P(n^{-1/3+\varepsilon}), \end{aligned}$$

where the order $O_P(\cdot)$ is uniform over $t \in [0, \tau]$. Notice that the integrands in the numerator are bounded conditionally independent with conditional mean zero when conditioning on $\mathcal{F}$. Then (A.7) can be shown by applying the empirical approximation analogous to the proof of (A.1).

Equation (A.8) follows from the definition of $g(\cdot)$, (A.6) and (A.7).

To show (A.9), recall that $h_0(t)$ is defined as $E[\tilde{w}(t)\exp\{\beta'Z(t)\}|Y\geq t]$. Use (A.1) and (A.7) and write

$$\frac{1}{n}\int_0^\tau \psi_n(t-s)\frac{\sum_{i=1}^n \tilde{w}_i(t)\,dN_i(t)}{\tilde{S}_0(t,\beta)}$$

$$= \frac{1}{n}\sum_{i=1}^n \int_0^\tau \frac{n^{1/3}\psi_n(t-s)\tilde{w}_i(t)\,dN_i(t)}{mf_n(t)h_0(t)} + O_P(n^{1/3+\varepsilon})$$

$$= \frac{1}{n}\sum_{i=1}^n \frac{n^{1/3}\psi_n(Y_i-s)\tilde{w}_i(Y_i)\delta_i}{mf_n(Y_i)h_0(Y_i)} + O_P(n^{1/3+\varepsilon})$$

(A.12)
$$= E\left[\frac{n^{1/3}\psi_n(Y-s)\tilde{w}(Y)\delta}{mf_n(Y)h_0(Y)}\right] + O_P(n^{-1/3+\varepsilon})$$

$$= \int_0^\tau \left[\frac{n^{1/3}\psi_n(t-s)\lambda_0(t)P(Y\geq t)}{mf_n(t)}\right] dt + O_P(n^{-1/3+\varepsilon})$$

$$= \frac{\lambda_0(s)P(Y\geq s)}{mf_n(s)} + O_P(n^{-1/3+\varepsilon}),$$

where $O_P(\cdot)$ is uniform over $s\in[0,\tau]$. In the above equations, the third equality can be proved analogous to the proof of (A.1). The details are omitted. That $\lambda_0(s)P(Y\geq s)/\{mf_n(s)\} = 1/\{mb(s)\}$ for $s\in(n^{-1/3},\tau-n^{-1/3})$ follows from the result of part 1 and the definitions of $b(\cdot)$ and $f_n(\cdot)$.

Equations (A.7)–(A.9) also hold if $\tilde{S}_k(t,\beta)$ is replaced by $S_k(t,\beta)$ because of (A.6). The proof of this lemma is complete. $\square$

PROOF OF THE THEOREM. Define $M(t) = N(t) - \int_0^t \exp\{\beta'Z(s)\}Y(s) \times \lambda_0(s)\,ds$. Let $M_i(t), i=1,2,\ldots,$ be the i.i.d. copies of $M(t)$,

$$g_n(t) = \frac{\sum_{i=1}^n \tilde{w}_i(t)Z_i(t)\exp\{\beta'Z_i(t)\}Y_i(t+)}{\sum_{i=1}^n \tilde{w}_i(t)\exp\{\beta'Z_i(t)\}Y_i(t+)}$$

and

$$\tilde{U}(\beta) = \sum_{i=1}^n \int_0^\tau \tilde{w}_i(t)\left\{Z_i(t) - \frac{\tilde{S}_1(t,\beta)}{\tilde{S}_0(t,\beta)}\right\} dN_i(t).$$

Then condition (vi) implies that, for any $\varepsilon_n \downarrow 0$,

(A.13) $$\sup_{|t-s|\leq\varepsilon_n, 0\leq t,s<\tau} |g_n(t) - g_n(s) - g(t) + g(s)| = o_P(n^{-1/2}).$$

The rest of the proof is divided into four steps.

*Step* 1 [To show $n^{-1/2}\tilde{U}(\beta) \to N(0, \Sigma)$]. Apply (A.1) and write

$$\tilde{U}(\beta) = \sum_{i=1}^{n} \int_0^\tau \tilde{w}_i(t)\{Z_i(t) - g_n(t-)\}\,dN_i(t)$$

$$- \sum_{i=1}^{n} \int_0^\tau \left\{\frac{\tilde{S}_1(t,\beta)}{\tilde{S}_0(t,\beta)} - g_n(t-)\right\}\sum_{i=1}^{n}\tilde{w}_i(t)\,dN_i(t)$$

$$= \sum_{i=1}^{n} \int_0^\tau \tilde{w}_i(t)\{Z_i(t) - g_n(t-)\}\,dM_i(t)$$

$$- \int_0^\tau \int_0^\tau \sum_{j \in R_s^*} \tilde{w}_j(s)\{Z_j(s) - g_n(t-)\}$$

$$\times \exp\{\beta'Z_j(s)\}\psi_n(t-s)\,d\bar{N}(s)\frac{\sum_{i=1}^n \tilde{w}_i(t)}{\tilde{S}_0(t,\beta)}\,dN_i(t)$$

$$= \sum_{i=1}^{n} \int_0^\tau \tilde{w}_i(t)[Z_i(t) - g_n(t-)]\,dM_i(t)$$

$$- \int_0^\tau \left[\int_0^\tau \sum_{j \in R_s^*} \tilde{w}_j(s)\{Z_j(s) - g_n(s)\}\right.$$

$$\left.\times \exp\{\beta'Z_j(s)\}\psi_n(t-s)\,d\bar{N}(s)\right]\frac{\sum_{i=1}^n \tilde{w}_i(t)\,dN_i(t)}{\tilde{S}_0(t,\beta)}$$

$$- \int_0^\tau \int_0^\tau \{g_n(s) - g_n(t-)\}\sum_{j \in R_s^*}\tilde{w}_j(s)\exp\{\beta'Z_j(s)\}$$

$$\times \left[\frac{\psi_n(t-s)}{\tilde{S}_0(t,\beta)}\sum_{i=1}^{n}\tilde{w}_i(t)\,dN_i(t)\right]d\bar{N}(s)$$

$$= \Xi_1 + \Xi_2 + \Xi_3, \quad \text{say.}$$

We first show $\Xi_3 = o_P(n^{1/2})$. In view of (A.13), it is seen that $\Xi_3$ differs by a term of order $o_P(n^{1/2})$ when $g_n(s) - g_n(t-)$ is replaced by $g(s) - g(t-)$. Notice that condition (v) ensures the differentiability of $g(\cdot)$. Therefore

$$\sup\{|g(s) - g(t)|\psi_n(t-s) : t, s \in [0, \tau]\} = O(n^{-1/3})$$

since $\psi(\cdot)$ has bounded support. Then, using the delta method, $\Xi_3$ can be reduced to

$$n^{4/3}\int_0^\tau \int_0^\tau \{g(s) - g(t)\}\psi_n(t-s)\,ds\,E[\tilde{w}(t)\exp\{\beta'Z(t)\}Y(t)]\lambda_0(t)\,dt$$

$$+ o_P(n^{1/2})$$

$$= O(n^{1/3}) + o_P(n^{1/2})$$
$$= o_P(n^{1/2})$$

by condition (v) and the result of part 1.

We next show the asymptotic normality of $\Xi_2$. Recall that $\mathcal{F}$ is the $\sigma$-algebra generated by $[\{Y_i, \delta_i, Z_i(\cdot)\} : i = 1, 2, \ldots]$. For every $t$, the integral in the brackets, conditioning on $\mathcal{F}$, has mean 0 and, when normalized by $\int_0^\tau \psi_n(t-s) \, d\bar{N}(s)$, can be shown to be $O_P(n^{-1/3+\varepsilon})$ uniformly over $t \in [0, \tau]$ for any fixed $\varepsilon > 0$. Thus, in view of (A.9), the delta method can be applied to show that $\Xi_2$ is

$$\int_0^\tau \sum_{j \in R_s^*} \tilde{w}_j(s)\{Z_j(s) - g_n(s)\}$$
$$\times \exp\{\beta' Z_j(s)\} \left[ \int_0^\tau \frac{\psi_n(t-s)}{\tilde{S}_0(t,\beta)} \sum_{i=1}^n \tilde{w}_i(t) \, dN_i(t) \right] d\bar{N}(s)$$
$$= \int_0^\tau \sum_{j \in R_s^*} \tilde{w}_j(s)\{Z_j(s) - g_n(s)\}$$
$$\times \exp\{\beta' Z_j(s)\} \frac{1}{mb(s)} \sum_{i=1}^n dN_i(s) + o_P(n^{1/2}),$$

where the main term, conditioning on $\mathcal{F}$, is the sum of bounded random variables with conditional mean 0. Therefore it converges at the rate $n^{1/2}$ to a normal distribution with mean 0 and asymptotic variance

$$\Sigma_0 \equiv \int_0^\tau E\left[\frac{\tilde{w}(t)^2}{mb(t)}\{Z(t) - g(t)\}^{\otimes 2} \exp\{2\beta' Z(t)\} Y(t)\right] \lambda_0(t) \, dt$$
$$= \int_0^\tau E[\tilde{w}(t)\{1 - \tilde{w}(t)\}\{Z(t) - g(t)\}^{\otimes 2} \exp\{\beta' Z(t)\} Y(t)] \lambda_0(t) \, dt.$$

Consider the first term $\Xi_1$ and write

$$\Xi_1 = \sum_{i=1}^n \int_0^\tau \tilde{w}_i(s)[Z_i(s) - g_n(s-)] \, dM_i(s)$$
$$= \sum_{i=1}^n \int_0^\tau \tilde{w}_i(s)[Z_i(s) - g(s)] \, dM_i(t) - \sum_{i=1}^n \int_0^\tau \tilde{w}_i(s)[g_n(s-) - g(s)] \, dM_i(t)$$
$$= \sum_{i=1}^n \int_0^\tau \tilde{w}_i(s)[Z_i(s) - g(s)] \, dM_i(t) + o_P(n^{1/2}).$$

The main term in the above expression is the sum of i.i.d. bounded random variables with mean 0 and variance $\Sigma_1 \equiv \mathrm{var}[\int_0^\tau \tilde{w}(t)\{Z(t) - g(t)\} \, dM(t)]$.

Therefore $\Xi_1$ converges at the rate $n^{1/2}$ to a normal distribution with mean 0 and variance $\Sigma_1$ defined above.

Last, combine the limits of the terms $\Xi_1, \Xi_2$ and $\Xi_3$, and notice that $\Xi_1$ is $\mathcal{F}$-measurable and that the asymptotic normality for $\Xi_2$ holds conditioning on $\mathcal{F}$. Observe that $\Sigma_1 = E[\tilde{w}(s)^2 \{Z(s) - g(s)\}^{\otimes 2} \exp\{\beta' Z(s)\} Y(s)] \lambda_0(s) \, ds$. It follows that

$$n^{-1/2}(\Xi_1 + \Xi_2 + \Xi_3) \to N(0, \Sigma_0 + \Sigma_1) = N(0, \Sigma).$$

*Step* 2 [To show $n^{-1/2} U(\beta) \to N(0, \Sigma)$]. The following equalities use the approximations (A.5)–(A.9) and the delta method:

$$\begin{aligned}
&U(\beta) - \tilde{U}(\beta) \\
&= \sum_{i=1}^n \int_0^\tau \left( \{w_i(t) - \tilde{w}_i(t)\} \left\{ Z_i(t) - \frac{\tilde{S}_1(t, \beta)}{\tilde{S}_0(t, \beta)} \right\} \right. \\
&\qquad\qquad \left. - w_i \left[ \frac{S_1(t, \beta)}{S_0(t, \beta)} - \frac{\tilde{S}_1(t, \beta)}{\tilde{S}_0(t, \beta)} \right] \right) dN_i(t) \\
&= \sum_{i=1}^n \int_0^\tau \left( \{w_i(t) - \tilde{w}_i(t)\} \{Z_i(t) - g(t)\} \right. \\
&\qquad\qquad - \tilde{w}_i(t) \left[ \frac{S_1(t, \beta) - \tilde{S}_1(t, \beta)}{\tilde{S}_0(t, \beta)} \right. \\
&\qquad\qquad \left. \left. - \frac{\tilde{S}_1(t, \beta)}{\tilde{S}_0^2(t, \beta)} \{S_0(t, \beta) - \tilde{S}_0(t, \beta)\} \right] \right) dN_i(t) + o_P(n^{1/2}) \\
&= \sum_{i=1}^n \int_0^\tau (\{w_i(t) - \tilde{w}_i(t)\}\{Z_i(t) - g(t)\}) \, dN_i(t) \\
&\quad - \int_0^\tau \frac{1}{\tilde{S}_0(t, \beta)} [S_1(t, \beta) - \tilde{S}_1(t, \beta) \\
&\qquad\qquad - g(t)\{S_0(t, \beta) - \tilde{S}_0(t, \beta)\}] \sum_{i=1}^n \tilde{w}_i(t) \, dN_i(t) + o_P(n^{1/2}) \\
&= \sum_{i=1}^n \int_0^\tau \{w_i(t) - \tilde{w}_i(t)\}\{Z_i(t) - g(t)\} \, dN_i(t) \\
&\quad - \int_0^\tau \left[ \sum_{j \in R_s^*} \{w_j(s) - \tilde{w}_j(s)\} \right. \\
&\qquad\qquad \times \{Z_j(s) - g(t)\} \exp\{\beta' Z_j(s)\} Y_j(s+)
\end{aligned}$$

$$\times \int_0^\tau \psi_n(t-s) \frac{\sum_{i=1}^n \tilde{w}_i(t)\, dN_i(t)}{\tilde{S}_0(t,\beta)} \Bigg] d\bar{N}(s) + o_P(n^{1/2})$$

$$= \sum_{i=1}^n \int_0^\tau \{w_i(t) - \tilde{w}_i(t)\}\{Z_i(t) - g(t)\}\, dM_i(t)$$

$$- \frac{1}{m} \int_0^\tau \Bigg[ \sum_{j=1}^n \Big\{ I(j \in R_s^*) - \frac{m}{n_s} \Big\}$$

$$\times \{w_j(s) - \tilde{w}_j(s)\}\{Z_j(s) - g(s)\}$$

$$\times \exp\{\beta' Z_j(s)\} Y_j(s+) \frac{n}{b(s)} \Bigg] d\bar{N}(s) + o_P(n^{1/2}).$$

The two main terms are also of order $o_P(n^{1/2})$, by observing the expressions (A.10) and (A.11) and an argument along the line of the proof of Theorem 4.1 of Sasieni ([1993a](#)). Hence $U(\beta) - \tilde{U}(\beta) = o_P(n^{1/2})$ and $n^{1/2} U(\beta) \to N(0, \Sigma)$ follows from Step 1.

*Step* 3 [To show the consistency of $\hat{\beta}$]. Let $B(\beta, \varepsilon_0)$ be the ball in $\mathbb{R}^d$ centered at $\beta$ with radius $\varepsilon_0 > 0$. In view of (A.6)–(A.8), one can show that, as $n \to \infty$ and then $\varepsilon_0 \to 0$,

(A.14) $$\sup_{v \in B(\beta, \varepsilon_0)} |-\dot{U}(v)/n - \Sigma| \to 0$$

in probability, where $v$ may be different for different elements of the matrix $\dot{U}(\cdot)$. Now, choose any small but fixed $\varepsilon_0 > 0$ and view $-U(\cdot)/n$ as a random mapping from $\mathbb{R}^d$ to $\mathbb{R}^d$. Then, since $\Sigma$ is assumed to be positive definite and thus invertible in condition (iv), (A.14) implies that, with probability tending to 1, the mapping $-U(\cdot)/n$ is a homeomorphism from $B(\beta, \varepsilon_0)$ to its image, denoted as $B_n$, which contains a ball of fixed radius. Since $U(\beta)/n = O_P(n^{-1/2})$ as proved in Step 2, $B_n$ contains $0 \in \mathbb{R}^d$ with probability tending to 1. Therefore, $\hat{\beta}$, as the unique solution of $U(\cdot) = 0$, is in the ball $B(\beta, \varepsilon_0)$ with probability tending to 1. The consistency of $\hat{\beta}$ is proved since $\varepsilon_0$ is arbitrary.

*Step* 4 [To show $n^{1/2}(\hat{\beta} - \beta) \to N(0, \Sigma^{-1})$]. It follows from the mean value theorem that

$$-U(\beta) = U(\hat{\beta}) - U(\beta) = \dot{U}(v)(\hat{\beta} - \beta),$$

where $v$ lies on the line segment joining $\hat{\beta}$ and $\beta$ but may be different for different elements of the matrix $\dot{U}(\cdot)$. Then the desired asymptotic normality of $\hat{\beta}$ follows from (A.14), the consistency of $\hat{\beta}$ proved in Step 3 and the asymptotic normality of $U(\beta)$ proved in Step 2. The proof is complete. $\square$

PROOF OF THE PROPOSITION. Define $\eta(t) = \{1 - \tilde{w}(t)\}\{Z(t) - g(t)\} + \{g(t) - \mu(t)\}/(m+1)$. Let $\eta_i(t), i = 1, \ldots, m+1$, be the i.i.d. copies of $\eta(t)$. The following calculations use the fact that $E[\{Z(t) - \mu(t)\} \exp\{\beta'Z(t)\}|Y \geq t] = E[\tilde{w}(t)\{Z(t) - g(t)\} \exp\{\beta'Z(t)\}|Y \geq t] = E[\{1 - \tilde{w}(t)\}\{Z(t) - g(t)\}|Y \geq t] = 0$. Write

$$\frac{1}{m+1} E\left(\sum_{k=1}^{m+1} \eta_k(t) \sum_{j=1}^{m+1} [\{Z_j'(t) - \mu'(t)\} \exp\{\beta'Z_j(t)\}] \Big| Y_1 \geq t, \ldots, Y_{m+1} \geq t\right)$$

$$= E[\eta(t)\{Z'(t) - \mu'(t)\} \exp\{\beta'Z(t)\}|Y \geq t]$$
$$= E[\{1 - \tilde{w}(t)\}\{Z(t) - g(t)\}\{Z'(t) - \mu'(t)\} \exp\{\beta'Z(t)\}|Y \geq t]$$
$$= E[\{Z(t) - \mu(t)\}^{\otimes 2} \exp\{\beta'Z(t)\}|Y \geq t]$$
$$\quad - E[\tilde{w}(t)\{Z(t) - g(t)\}^{\otimes 2} \exp\{\beta'Z(t)\}|Y \geq t]$$
$$= H_1(t) - H_2(t), \quad \text{say.}$$

Similarly,

$$\frac{1}{m+1} E\left[\left\{\sum_{k=1}^{m+1} \eta_k(t)\right\}^{\otimes 2} \sum_{j=1}^{m+1} \exp\{\beta'Z_j(t)\} \Big| Y_1 \geq t, \ldots, Y_{m+1} \geq t\right]$$

$$= E[\eta(t)^{\otimes 2} \exp\{\beta'Z(t)\}|Y \geq t] + mE\{\eta(t)^{\otimes 2}|Y \geq t\}b(t)$$
$$\quad + 2mE[\eta(t) \exp\{\beta'Z(t)\}|Y \geq t]E\{\eta'(t)|Y \geq t\}$$
$$\quad + m(m-1)[E\{\eta(t)|Y \geq t\}]^{\otimes 2} b(t)$$
$$= E(\eta(t)^{\otimes 2}[\exp\{\beta'Z(t)\} + mb(t)]|Y \geq t)$$
$$\quad + \frac{m}{m+1}\left\{2\left(-1 + \frac{1}{m+1}\right) + \frac{m-1}{m+1}\right\}\{g(t) - \mu(t)\}^{\otimes 2} b(t)$$
$$= E\left[\frac{\eta(t)^{\otimes 2} \exp\{\beta'Z(t)\}}{1 - \tilde{w}(t)}\Big| Y \geq t\right] - \frac{m}{m+1}\{g(t) - \mu(t)\}^{\otimes 2} b(t)$$
$$= E[\{1 - \tilde{w}(t)\}\{Z(t) - g(t)\}^{\otimes 2} \exp\{\beta'Z(t)\}|Y \geq t]$$
$$\quad - \{g(t) - \mu(t)\}^{\otimes 2} b(t)$$
$$= E[\{Z(t) - \mu(t)\}^{\otimes 2} \exp\{\beta'Z(t)\}|Y \geq t]$$
$$\quad - E[\tilde{w}(t)\{Z(t) - g(t)\}^{\otimes 2} \exp\{\beta'Z(t)\}|Y \geq t]$$
$$= H_1(t) - H_2(t).$$

Then it follows from (the matrix version of) the Cauchy–Schwarz inequality that

$$H_3(t) \equiv \frac{1}{m+1} E\left( \frac{[\sum_{j=1}^{m+1}\{Z_j(t) - \mu(t)\} \exp\{\beta' Z_j(t)\}]^{\otimes 2}}{\sum_{j=1}^{m+1} \exp\{\beta' Z_j(t)\}} \bigg| Y_1 \geq t, \ldots, Y_{m+1} \geq t \right)$$
$$\geq H_1(t) - H_2(t),$$

and equality holds if and only if

$$\mathrm{pr}\left[ \sum_{j=1}^{m+1} \eta_j(t) = h(t) \frac{\sum_{j=1}^{m+1}\{Z_j(t) - \mu(t)\}\exp\{\beta' Z_j(t)\}}{\sum_{j=1}^{m+1} \exp\{\beta' Z_j(t)\}} \bigg| Y_1 \geq t, \ldots, Y_{m+1} \geq t \right] = 1,$$

where $h(t)$ is a nonrandom function. This equality holds only when the conditional distribution of $\sum_{j=1}^{m+1} \exp\{\beta' Z_j(t)\}$ given $Y_1 \geq t, \ldots, Y_{m+1} \geq t$ is degenerate. If the above equality holds for all $t \in [0, \tau]$ except for a Lebesgue measure 0 set, then $\beta = 0$. Observe that $\Sigma_C = \int_0^\tau H_1(t) \mathrm{pr}(Y \geq t) \lambda_0(t) \, dt$, $\Sigma = \int_0^\tau H_2(t) \mathrm{pr}(Y \geq t) \lambda_0(t) \, dt$ and $\Sigma_a = \int_0^\tau H_3(t) \mathrm{pr}(Y \geq t) \lambda_0(t) \, dt$. Then it follows that $\Sigma_a \geq \Sigma_C - \Sigma$, or, equivalently, $\Sigma_P^{-1} = (\Sigma_C - \Sigma_a)^{-1} \geq \Sigma^{-1}$. And $\Sigma_P^{-1} = \Sigma^{-1}$ if and only if $\beta = 0$. The proof is complete. $\square$

**Acknowledgments.** The author is deeply grateful to Professor Langholz for his enlightening comments and discussions and for his generosity in sharing his many manuscripts, and to a referee and Ms. Joanne Zhong for their careful reading of this paper and detailed suggestions.

## REFERENCES

ANDERSEN, P. K. and GILL, R. D. (1982). Cox's regression model for counting processes: A large sample study. *Ann. Statist.* **10** 1100–1120. MR673646

BORGAN, Ø., GOLDSTEIN, L. and LANGHOLZ, B. (1995). Methods for the analysis of sampled cohort data in the Cox proportional hazards model. *Ann. Statist.* **23** 1749–1778. MR1370306

BORGAN, Ø. and OLSEN, E. F. (1999). The efficiency of simple and counter-matched nested case-control sampling. *Scand. J. Statist.* **26** 493–509. MR1734258

BRESLOW, N. E. (1996). Statistics in epidemiology: The case-control study. *J. Amer. Statist. Assoc.* **91** 14–28. MR1394064

CHEN, K. (2001). Generalized case-cohort sampling. *J. R. Stat. Soc. Ser. B Stat. Methodol.* **63** 791–809. MR1872067

CHEN, K. and LO, S.-H. (1999). Case-cohort and case-control analysis with Cox's model. *Biometrika* **86** 755–764. MR1741975

GOLDSTEIN, L. and LANGHOLZ, B. (1992). Asymptotic theory for nested case-control sampling in the Cox regression model. *Ann. Statist.* **20** 1903–1928. MR1193318

LANGHOLZ, B. and GOLDSTEIN, L. (1996). Risk set sampling in epidemiologic cohort studies. *Statist. Sci.* **11** 35–53.

Langholz, B. and Thomas, D. C. (1990). Nested case-control and case-cohort methods of sampling from a cohort: A critical comparison. *Amer. J. Epidemiology* **131** 169–176.

Langholz, B. and Thomas, D. C. (1991). Efficiency of cohort sampling designs: Some surprising results. *Biometrics* **47** 1563–1571.

Oakes, D. (1981). Survival times: Aspects of partial likelihood (with discussion). *Internat. Statist. Review* **49** 235–264. MR651473

Pollard, D. (1990). *Empirical Processes*: *Theory and Applications*. IMS, Hayward, CA. MR1089429

Robins, J. M., Rotnitzky, A. and Zhao, L. P. (1994). Estimation of regression coefficients when some regressors are not always observed. *J. Amer. Statist. Assoc.* **89** 846–866. MR1294730

Samuelsen, S. O. (1997). A pseudo-likelihood approach to analysis of nested case-control studies. *Biometrika* **84** 379–394. MR1467054

Sasieni, P. (1993a). Some new estimators for Cox regression. *Ann. Statist.* **21** 1721–1759. MR1245766

Sasieni, P. (1993b). Maximum weighted partial likelihood estimators for the Cox model. *J. Amer. Statist. Assoc.* **88** 144–152.

Suissa, S., Edwardes, M. and Biovin, J.-F. (1998). External comparisons from nested case-control designs. *Epidemiology* **9** 72–78.

Thomas, D. C. (1977). Addendum to "Methods of cohort analysis: Appraisal by application to asbestos mining," by F. D. K. Liddell, J. C. McDonald and D. C. Thomas. *J. Roy. Statist. Soc. Ser. A* **140** 483–485.

van der Vaart, A. W. and Wellner, J. A. (1996). *Weak Convergence and Empirical Processes*: *With Applications to Statistics*. Springer, New York. MR1385671

Department of Mathematics
HKUST
Clear Water Bay
Kowloon
Hong Kong
e-mail: makchen@ust.hk